  \def\addforarchive{\begin{picture}(0,0)
                     \put(220,-80){
                     \put(60,250){{\small\sf math.CT/0602079}}
                     \put(60,238){{\small\sf KCL-MTH-06-01}}
                     \put(60,227){{\small\sf
                     Hamburger$\;$Beitr.$\;$zur$\;$Math.$\;$Nr.$\;$230}}
                     \put(60,215){{\small\sf ZMP-HH/05-29}}}
                     \end{picture}
                     \vspace*{-2em}
                     }
  \def\Bi{\BiA} \def\BiSize{FfFFF7}

\documentclass{icmart}
\contact[schweigert@math.uni-hamburg.de]
        {Organisationseinheit Mathematik, Universit\"at Hamburg,
        Bundesstra{\rm\ss{}}e 55, D\,--\,20146 Hamburg}
\contact[jfuchs@fuchs.tekn.kau.se]
        {Avdelning fysik, Karlstads Universitet,
        Universitetsgatan~5, S\,--\,65188 Karlstad}
\contact[ingo@mth.kcl.ac.uk]
        {Department of Mathematics, King's College London, Strand,
        London WC2R 2LS}

\theoremstyle{definition}

\newcommand{\Hom}{\operatorname{Hom}}

\def\bysame         {------}

\def\calc          {\ensuremath{\mathcal C}}
\def\cald          {\ensuremath{\mathcal D}}
\def\calg          {\ensuremath{\mathcal G}}

\def\call          {\ensuremath{\mathcal L}}
\def\calm          {\ensuremath{\mathcal M}}

\def\calv          {\ensuremath{\mathcal V}}
\def\Bimod         {\mbox{\sl Bimod}}
\def\C             {\ensuremath{\mathbb C}}
\def\Cardy         {\mbox{\small\sl I}}
\def\CA            {\ensuremath{\mathcal C_{\!A}}}
\def\cft           {conformal field theory}
\def\cfts          {conformal field theories}

\def\cir           {\,{\circ}\,}
\def\Cf            {\mbox{\sl Cor}}
\def\eps           {\varepsilon}
\def\eq            {\,{=}\,}
\def\FROB          {\ensuremath{\mathcal F}\!\mbox{\sl rob}}
\def\Hom           {{\rm Hom}}
\def\HOM           {\ensuremath{\mathcal H}\mbox{\sl om}}
\def\hatX          {\ensuremath{\widehat{\mathrm X}}}
\def\hy            {$\mbox{-\hspace{-.66 mm}-}$\linebreak[0]}
\def\id            {\mbox{\sl id}}
\def\iN            {\,{\in}\,}
\def\io            {\mbox{\sl\i}}
\def\M             {\ensuremath{\mathrm M}}
\def\Mod           {\mbox{\sl Mod}}
\def\N             {\ensuremath{\mathbb N}}
\def\one           {{\bf1}}
\def\oti           {\,{\otimes}\,}

\def\qft           {\mbox{\sl qft}}
\def\R             {\ensuremath{\mathbb R}}

\def\tft           {\mbox{\sl tft}}
\def\X             {\ensuremath{\mathrm X}}
\def\Y             {\ensuremath{\mathrm Y}}
\def\Z             {\ensuremath{\mathbb Z}}

\newcommand\void[1]{}

\title[Categorification and correlation functions in conformal field theory]
      {Categorification and correlation functions in conformal field theory}
\author[C.\ Schweigert, J.\ Fuchs and I.\ Runkel]
       {Christoph Schweigert, J\"urgen Fuchs, and Ingo Runkel 
\thanks{~%
       C.S.\ is supported by the DFG project SCHW 1162/1-1,
       and J.F.\ by VR under project no.\ 621--2003--2385.
} }

\begin{document}

\addforarchive

\begin{abstract}
A modular tensor category provides the appropriate data for the construction 
of a three-dimensional topological field theory. We describe the following 
analogue for two-dimensional conformal field theories: a 2-category whose 
objects are symmetric special Frobenius algebras in a modular tensor category 
and whose morphisms are categories of bimodules. 
This 2-category provides sufficient ingredients for constructing all
correlation functions of a two-dimensional rational conformal field theory.
The bimodules have the physical interpretation of chiral data, boundary 
conditions, and topological defect lines of this theory.
\end{abstract} 



\maketitle


\section{Quantum field theories as functors}

In approaches to quantum field theory that are based on the concepts
of fields and states, the utility of categories and functors is by
now well-established. The following pattern has been recognized: There is a 
geometric category $\calg$ which, for every concrete model,
must be suitably ``decorated''. The decoration is achieved with the help of 
objects and morphisms from another category $\calc$. For known classes of 
quantum field theories, the decoration category  $\calc$ typically 
has a representation-theoretic origin -- the reader is encouraged to think of 
it as the representation category of some algebraic object, like a quantum 
group, a loop group, a vertex algebra, a net of observable algebras, etc. 
This way one obtains a decorated geometric category $\calg_\calc$. The quantum 
field theory can then be formulated as a (tensor) functor $\qft_{\calc}$ from 
$\calg_\calc$ to some category of vector spaces.
In this contribution, we mainly consider cases for which this latter category 
is the tensor category of finite-dimensional complex vector spaces. 

A prototypical example for this pattern is provided by topological quantum field 
theories (TFTs). For such theories, the geometric category $\calg$ is based on 
a cobordism category: its objects are $d{-}1$-dimensional topological manifolds 
without boundary. It is convenient to include two types of morphisms
\cite{stte}: homeomorphisms of $d{-}1$-dimensional manifolds, and cobordisms. 
A cobordism $\M{:}\ \Y_{\!1}\,{\to}\,\Y_{\!2}$ is a $d$-di\-men\-si\-o\-nal 
topological manifold $\M$ together with a parametrization of its boundary given 
by a homeomorphism $\partial\M \,{\stackrel\cong\to}\, \overline{\Y_{\!1}} 
\,{\sqcup}\, \Y_{\!2}$, where $\overline{\Y_{\!1}}$ has the orientation opposite
to the one of $\Y_{\!1}$. The composition of morphisms is by concatenation, by 
gluing, and by changing the parametrization of the boundary, respectively.
Cobordisms that coincide by a homeomorphism of the $d$-manifold $\M$ 
restricting to the identity on $\partial\M$ must be identified.

In the simplest case, a topological field theory thus associates to a closed 
$d{-}1$-di\-mensional manifold $\X$ a vector space $\qft_\calc(\X)$, and to a 
homeomorphism or a co\-bordism $\M{:}\ \Y_{\!1}\,{\to}\,\Y_{\!2}$ a linear map
$$
\qft_\calc(\M):\quad \qft_\calc(\Y_{\!1}) \to \qft_\calc(\Y_{\!2}) \,.
$$
The assignment $\qft_\calc$ is required to be a (strict) tensor functor.
This requirement implies the usual axioms (cf.\ e.g.\ \cite{TUra})
of naturality, multiplicativity, functoriality and normalization.

There also exists a path-integral approach to certain classes of topological 
field theories. Its relation to the categorical framework described above is as 
follows: One can think of the vector space $\qft_\calc(\partial_-\M)$ as the
space of (equivalence classes of) possible initial data for ``fields'' in the 
path integral, and of $\qft_\calc(\partial_+\M)$ as the possible final data. 
The matrix elements of the linear map $\qft_\calc(\M)$ are then the transition 
amplitudes for fixed initial and final values of the fields. 

This picture is still oversimplified. In particular, it turns out that it is 
natural to enrich also the geometric category $\calg$ over the category of 
complex vector spaces. As a consequence, when studying functors on $\calg$,
one should then consider also projective functors.
These issues, which are closely related to anomalies in quantum field theory 
will, however, be suppressed in this note.

A prominent class of examples of $d\eq3\,$-dimensional topological field 
theories arises from Chern-Simons field theories. 
For $G$ a simple connected and simply-connected complex Lie group, consider
holomorphic $G$-bundles on a closed two-manifold \X\ of genus $g$ with complex 
structure. Pick a generator $\call$ for the Picard group of the moduli space 
$\calm^G_{\X}$ of such bundles. 
For every fixed $k\iN\N$, upon changing the complex structure of \X, the 
vector spaces $H^0(\calm^G_{\X},\call^{\otimes k})$ fit together into a vector 
bundle with projectively flat connection over the moduli space $\calm^G_{g}$ of 
curves of genus $g$. The complex modular functor \cite{BAki} associates 
these bundles to \X; these bundles and their monodromies provide a formalization
of all aspects of the chiral level-$k$ Wess-Zumino-Witten (WZW) conformal field 
theory for $G$ that are needed for the discussions in the subsequent sections.

As a next step, it is natural to extend the formalism by allowing for marked 
points with additional structure on the two-manifolds. {}From a field 
theoretical point of view, this is motivated by the desire to account for 
insertions of fields. In the case of Chern-Simons theory, the additional 
structure amounts to specifying parabolic structures at the marked points. 
The marked points thus have to carry labels, 
which we will identify in a moment as objects of a decoration category $\calc$.  

This structure must be extended to the geometric morphisms: Maps of 
$2$-di\-men\-sional manifolds are required to preserve marked points and the 
decoration in $\calc$. The decoration of the $2$-di\-mensional manifolds 
is extended to the $3$-di\-men\-sional manifolds \M\ underlying cobordisms 
by supplying them with oriented (ribbon) graphs ending on (arcs through) the 
marked points on $\partial_\pm\M$. The ribbon graph is allowed to have 
vertices with a finite number of ingoing and outgoing ribbons. 
{}From the construction of invariants of knots and links, it is
known that this enforces a categorification of the set of labels:
$\calc$ must be a ribbon category, i.e.\ a braided sovereign tensor category.
In particular, the vertices of the graph are to be labeled by morphisms in the 
decoration category $\calc$.

This appraoch has been very fruitful and
has, in particular, made a rigorous construction of Chern-Simons theory
possible \cite{retu,TUra}. The extension from invariants of links in $\R^3$ to
link invariants in arbitrary oriented three-ma\-nifolds has revealed an
important subclass of tensor categories: modular tensor categories. 

For the purposes of the present contribution, we adopt the following definition 
of a modular tensor category: it is an abelian, $\C$-linear, semi-simple ribbon
category with a finite number of isomorphism classes of simple objects. The 
tensor unit $\one$ is required to be simple, and the braiding must be 
nondegenerate in the sense that the natural transformations of the identity 
functor on $\calc$ are controlled by the fusion ring $K_0(\calc)$:
$$
{\rm End}(\mbox{\sl Id}_{\calc}) \cong  K_0(\calc)\otimes_{\Z}^{} \C \,. 
$$

The relation between Chern-Simons theory and chiral Wess-Zumino-Witten theory
\cite{witt27} was a first indication that modular tensor categories also
constitute the appropriate mathematical formalization of the chiral data 
\cite{mose3,frki} of a conformal field theory. Recent progress in 
representation theory has made this idea much more precise; 
for the following classes of representation categories it has been 
established that they carry the structure of a modular tensor category:
   \def\leftmargini{1.1em}~\\[-1.4em]
   \begin{itemize}
\item
The representation category of a connected ribbon factorizable
weak Hopf algebra over $\C$ (or, more generally, over an algebraically
closed field $\Bbbk$) with a Haar integral \cite{nitv}.
\item
The category of unitary representations of the double of a connected $C^*$ 
weak Hopf algebra \cite{nitv}.
\item
The category of local sectors of a finite-index net of von Neumann algebras on 
the real line, if the net is strongly additive (which for conformal nets is 
equivalent to Haag duality) and has the split property \cite{kalm}.
\\
(In this example and in the previous one, one 
obtains unitary modular tensor categories.) 
\item
The representation category of a self-dual vertex algebra that obeys Zhu's 
$C_2$ cofiniteness condition and certain conditions on its homogeneous 
subspaces, provided that this category is semisimple \cite{huan21}.
\end{itemize}
The last two entries in this list correspond to two different mathematical 
formalizations of chiral conformal field theories. The results of \cite{kalm} 
and \cite{huan21} therefore justify the point of view that modular tensor
categories constitute an \oe{}cumenic formalization of the chiral data of a 
conformal field theory.


\section{Two-dimensional conformal field theories}

Three-dimensional topological field theory will indeed appear as a tool in 
constructions below. Our main interest here is, however, in a different class 
of quantum field theories: full (and in particular local) two-dimensional 
conformal field theories, or CFTs, for short.

For these theories, the geometric category of interest is the category of 
two-di\-mensional conformal manifolds, possibly with non-empty boundary. This 
fact already indicates that full conformal field theories are different 
from the chiral conformal field theories that we encountered in the last section 
for the case of WZW theories. Morphisms in the category of conformal manifolds
are maps respecting the conformal structure. Actually, 
there are two different types of full conformal field theories, corresponding 
to two different geometric categories: One considers either oriented conformal 
manifolds, leading to a category $\calg^{{\rm or}}$, or unoriented
manifolds, leading to a different geometric category $\calg^{{\rm unor}}$. 
As morphisms, we admit maps that preserve the respective structure.
(In the application of conformal field theory to string theory, 
$\calg^{{\rm or}}$ plays a role in superstrings of ``type II'', while 
$\calg^{{\rm unor}}$ appears in superstring theories of ``type I''.)

As in the case of topological field theories, the geometric category needs to 
be decorated. To find the appropriate decoration data, we first discuss 
a physical structure that is known to be present in specific classes of models 
and that our approach to conformal field theory should take into account:
   \def\leftmargini{1.1em}~\\[-1.2em]
   \begin{itemize}\addtolength\itemsep{2pt}
\item Whenever a two-manifold $\X$ has a boundary, one expects that it is
necessary to specify boundary conditions. In a path integral approach, a boundary
condition is a prescription for the boundary values of fields that appear in 
the Lagrangian. Here, a more abstract approach is adequate: we take the possible
boundary conditions to be the objects of a decoration category $\calm$. This 
constitutes again a categorification of the decoration data. It can be motivated 
further by the observation that insertions of boundary fields can change the 
boundary condition; they will be related to morphisms of the category $\calm$.
\end{itemize}

\noindent
A second structure, which in the literature has received much less attention
than boundary conditions, turns out to provide crucial clues for the 
construction of full conformal field theories:
   \def\leftmargini{1.1em}~\\[-1.2em]
   \begin{itemize}\addtolength\itemsep{2pt}
\item Conformal field theories can have topological defect lines. 
Such defect lines have e.g.\ been known for the Ising model for a long time: 
This CFT describes the continuum limit of a lattice model with $\Z_2$-valued
variables at the vertices of a two-dimensional lattice and with ferromagnetic 
couplings along its edges. A defect line is obtained when one changes the 
couplings on all edges that 
intersect a given line in the lattice from ferromagnetic to antiferromagnetic.

In the continuum limit, such a defect line can be described by a condition on 
the values of bulk fields at the defect line. In particular, when crossing a 
defect line, the correlation function of a bulk field can acquire a branch cut. 
Indeed, the defect lines we are interested in behave very much like branch 
cuts: they are topological in the sense that their precise location does not 
matter. In a field theoretic language, this is attributed to the fact that the 
stress-energy tensor of the theory is required to be smooth across defect lines.
As in the case of boundary conditions, in our framework it is not desirable to 
express defect lines through conditions on the values of fields. Instead, we 
anticipate again a categorification of the decoration data and label the 
possible types of defect lines by objects in yet another decoration category 
$\cald$.
\end{itemize}

There is a natural notion of fusion of defect lines, see e.g.\ \cite{pezu5,ffrs3}. 
Accordingly, $\cald$ will be a tensor category.  Also, 
to take into account the topological nature of defect lines, we assume that
the tensor category $\cald$ has dualities and that it is even sovereign. 
In contrast, there is no natural notion of a braiding of defect lines, so 
$\cald$ is, in general, not a ribbon category.

The two decoration categories -- $\cald$ for the defects and $\calm$ for the
boundary conditions -- are themselves related. One can fuse
a defect line to a boundary condition, thereby
obtaining another boundary condition; see e.g.\ \cite{grwa}.
This endows the category $\calm$ of boundary conditions with the
structure of a module category over $\cald$, i.e.\ there is a bifunctor
$$ 
\otimes :\quad \cald \times \calm \,\to\, \calm 
$$
which has (mixed) associativity constraints obeying a mixed pentagon identity.

The structure just unraveled -- a tensor category $\cald$ together with a 
module category $\calm$ over $\cald$ -- calls for the following natural 
extension: One should also consider the category of module functors, i.e.\ the 
category $\calc$ whose objects are endofunctors of $\calm$ that are compatible 
with the structure of a module category over $\cald$ and whose morphisms are 
natural transformations between such functors. The concatenation of functors 
naturally endows $\calc$ with a product so that $\calc$ is a tensor category. 

A recent insight is the following: In the application to 
two-dimensional conformal field theory, the category $\calc$ obtained this way
is equivalent to the category of chiral data that we discussed in Section 1\,!
There is a particularly amenable subclass of conformal field theories, called
{\em rational\/} conformal field theories ({\em R\/}CFTs), which can be 
rigorously discussed on the basis of this idea. For these theories, the 
category $\calc$ of chiral data has the structure of a modular tensor category.
In this case the idea can be exploited to arrive at a construction of
correlation functions (see Section \ref{CORF}) of rational conformal
field theories that is based on three-dimensional topological field
theory. This TFT approach to RCFT correlators will be presented 
in Section \ref{TFTRCFT} below.

\section{The 2-category of Frobenius algebras}

In practice, one frequently takes an opposite point of view: Instead of 
obtaining $\calc$ as a functor category, one starts from
some knowledge about the chiral symmetries of a conformal field theory. This 
allows one to use the representation-theoretic results mentioned in section 1 
so as to get a modular tensor category $\calc$ describing the chiral data
of the theory. Afterwards one realizes that the category $\calm$ of boundary
conditions is also a (right-) module category over $\calc$. General arguments 
involving internal Hom's \cite{ostr} together with specific properties relevant 
to rational conformal field theories then imply that in the tensor category 
$\calc$ there exists an associative algebra $A$ such that $\calm$ is equivalent
to the category \CA\ of left $A$-modules in $\calc$. By similar arguments one 
concludes that the category $\cald$ is equivalent to the category of 
$A$-bimodules. Additional constraints, in particular the non-degeneracy of the 
two-point functions of boundary fields on a disk, lead to further conditions 
on this algebra \cite{fuRs4}: $A$ must be a symmetric special Frobenius 
algebra. Owing to these insights we are able to use a generalization\,%
  \footnote{~When $\calc$ is the modular tensor category of finite-dimensional
  complex vector spaces, the CFT is a topological CFT. In particular,
  $A$ is then an ordinary Frobenius algebra. This case has served as a toy
  model for conformal field theories, see e.g.\ \cite{sega13,moor13}.}
of the theory of Frobenius algebras to braided tensor categories as a powerful 
tool to analyze (rational) conformal field theory. The algebraic theory in the 
braided setting is, however, genuinely richer; see \cite{ffrs} for a 
discussion of some new phenomena.

A {\em Frobenius\/} algebra $A\eq(A,m,\eta,\Delta,\eps)$ in $\calc$ is, by 
definition, an object of $\calc$ carrying the structures of a unital associative
algebra $(A,m,\eta)$ and of a counital coassociative coalgebra $(A,\Delta,\eps)$
in $\calc$, with the algebra and coalgebra structures satisfying the 
compatibility requirement that the coproduct $\Delta{:}\ A \,{\to}\, A\oti A$ 
is a morphism of $A$-bimodules (or, equivalently, that the product $m{:}\ A\oti 
A\,{\to}\,A$ is a morphism of $A$-bi-comodules). A Frobenius algebra is called
{\em special\/} iff the coproduct is a right-in\-ver\-se to the product --
this means in particular that the algebra is separable -- and a nonvanishing 
multiple of the unit $\eta{:}\ \one\,{\to}\,A$ is a right-inverse to the counit
$\eps{:}\ A\,{\to}\,\one$. There are two isomorphisms $A\,{\to}\, A^\vee$ that 
are naturally induced by product, counit and duality; $A$ is called 
{\em symmetric\/} iff these two isomorphisms coincide.  

Two algebras in a tensor category $\calc$ are called Morita equivalent iff their
representation categories are equivalent as module categories over $\calc$. 
Since the algebra $A$ is characterized by the requirement that \CA\ is 
equivalent to the given decoration category $\calm$, it is clear that only the 
Morita class of the algebra should matter. It is a non-trivial internal 
consistency check on the constructions to be presented in Section
\ref{TFTRCFT} that this is indeed the case.

Taking the modular tensor category $\calc$ as the starting point,
the following further generalization of the setup (compare also 
\cite{muge8,yama12,laud}) is now natural:\,%
  \footnote{~We are grateful to Urs Schreiber for discussions on this point.}
Consider the set of {\em all\/} (symmetric special) Frobenius algebras in 
$\calc$. This gives rise to a {\em family\/} of full conformal field 
theories that are based on the same chiral data. And again we categorify the 
structure: we introduce a 2-category $\FROB_\calc$ whose objects are symmetric 
special Frobenius algebras in $\calc$. The morphisms of a 2-category are 
categories; $\HOM(A,A')$ is simply the category of $A$-$A'$-bimodules. 
The 2-category $\FROB_\calc$ has a distinguished object $\Cardy\,$: as an 
object of $\calc$, $\Cardy$ is just the tensor unit, which is naturally a 
symmetric special Frobenius algebra. Because of the considerations 
in \cite{card9}, the full conformal field theory corresponding to $\Cardy$ is
often referred to as the ``Cardy case''; for this case a construction of the 
correlators in the spirit of Section \ref{CORF} was given in \cite{fffs3}.

We are now in a position to attribute a physical interpretation to the
morphisms of $\FROB_\calc$. $\HOM(\Cardy,\Cardy)$ is naturally
identified with the tensor category $\calc$ of chiral data. Further, for any
$A$ the tensor category $\HOM(A,A)$ describes topological defects in the full
conformal field theory associated to $A$; more generally, the category 
$\HOM(A,A')$ accounts for topological defect lines that separate two different
conformal field theories which share the same chiral data. Finally, the 
category $\HOM(\Cardy,A)$ also describes boundary conditions for the full
conformal field theory labeled by $A$.

We have thus learned that the decoration data of a family of full 
rational conformal field theories based on the same chiral data are described 
by a 2-category. This nicely fits with insight gained in other contexts:
\begin{itemize}
\item 
2-categories appear in recent approaches to elliptic objects \cite{badR,stte}.
\item 
Hermitian bundle gerbes, which appear naturally in a semi-classical description
of WZW conformal field theories \cite{gaRe}, form a 2-category \cite{steV}.
\end{itemize}
Unfortunately, at the time of writing, a unified approach to conformal field
theories based on 2-categories has not been established yet. For this reason, 
in the sequel we will not be able to use this language systematically.

\medskip

We close this section with two further comments. First, so far we have 
discussed the decoration data relevant to the oriented geometric category
$\calg^{{\rm or}}$. For the unoriented geometric category $\calg^{{\rm unor}}$,
additional structure on the relevant Frobenius algebra is needed:
$A$ must then be a {\em Jandl\/} algebra, that is, a symmetric special Frobenius
algebra coming with an algebra isomorphism $A \,{\stackrel\cong\to}\, 
A^{{\rm opp}}$ that squares to the twist on $A$. This turns out to be the 
appropriate generalization of the notion of an algebra with involution to 
braided tensor categories. We refrain from discussing this issue in the 
present contribution, but rather refer to \cite{fuRs8,fuRs11} for details.

Second, the general situation encountered above -- a module category $\calm$ 
over a tensor category $\calc$ -- naturally appears in various other contexts
as well:  
   \def\leftmargini{1.1em}~\\[-1.2em]
   \begin{itemize}\addtolength\itemsep{2pt}
\item 
The left modules over a weak Hopf algebra $H$ form a tensor
category $\calc\,{=}\,H$-\Mod. In a weak Hopf algebra, one can identify two
subalgebras $H_s$ and $H_t$ that are each other's opposed algebras \cite{bons}.
Forgetful functors thus endow any $H$-module with the structure of an 
$H_t$-bimodule; one even obtains a tensor functor from $H$-\Mod\ to 
$H_t$-\Bimod. The usual tensor product of an $H_t$-bimodule and an
$H_t$-left module endows the category of $H_t$-modules with the structure
of a module category over $H_t$-\Bimod\ and thus over $H$-\Mod. 

Weak Hopf algebras have indeed been proposed \cite{pezu6} as a framework for 
rational conformal field theories. Unfortunately, such a description must cope 
with two problems: First, to account for a braiding on $\calc$, one must work 
with an $R$-matrix on $H$; not too surprisingly, this is technically involved, 
and indeed not very much is known about $R$-matrices on weak Hopf algebras.
Secondly, given a tensor category $\calc$ (describing the chiral data), there 
does not exist a canonical Hopf algebra such that $H$-\Mod\ is equivalent to
$\calc$. Rather, as an additional datum, a fiber functor to $H_t$-bimodules 
needs to be chosen. A physical interpretation of this datum is unclear. On the 
other hand, Hopf algebras are still useful in the analysis of full rational CFT:
Their Hochschild cohomology was used \cite{etno} to compute the Davy\-dov-Yetter 
cohomology of the pair $(\calc,\calm)$; from the vanishing of this cohomology, 
rigidity properties of rational conformal field theories follow.
\item 
Weak Hopf algebras also appear in the study of inclusions of subfactors. For 
a review and further references, we refer to Sections 8 and 9 of \cite{niva2}.
\item 
The same category-theoretic structures have been recovered in the theory of 
vertex algebras from so-called open-string vertex algebras \cite{huko} which 
are, in particular, extensions of ordinary vertex algebras.
\end{itemize}
Not surprisingly, some of the structures that will be encountered
in the rest of this paper also have counterparts in the context of weak
Hopf algebras, of nets of subfactors, or of open-string vertex algebras.


\section{Correlation functions}\label{CORF}

The observations made in the preceding section raises the question whether 
one can construct a full rational conformal field theory by using a
modular tensor category $\calc$ and the 2-category $\FROB_\calc$ as an input.
These data should then in particular encode 
information about the correlation functions of the conformal field theory.

To decide this question, it is helpful to reformulate first the geometric 
categories $\calg^{{\rm or}}$ and $\calg^{{\rm unor}}$. This is achieved with 
the help of a crucial aspect of complex geometry in {\em two\/} dimensions: 
a complex structure on a two-dimensional manifold is equivalent to a conformal 
structure and the choice of an orientation. The complex double \hatX\ of a 
conformal manifold \X\ is a two-sheeted cover of \X\ whose points are pairs 
consisting of a point $p\iN\X$ and a local orientation at $p$. In view of 
the previous comment, it is clear that \hatX\ is a complex curve. We have thus 
associated to any object \X\ of the geometric categories $\calg^{{\rm unor}}$ 
and $\calg^{{\rm or}}$ a complex curve \hatX\ that comes with an
orientation reversing involution $\sigma$ such that the quotient $\hatX/\langle 
\sigma\rangle$ is naturally isomorphic to \X, and we have a canonical projection
$$  
\pi: \ \ \hatX \,\mapsto\,\X \cong \hatX\,/ \langle \sigma\rangle \,. 
$$
The set of fixed points of $\sigma$ is just the preimage  under $\pi$ of the 
boundary $\partial\X$.

To be able to use the tools of complex geometry, we therefore reformulate our
geometric categories as follows: in the case of $\calg^{{\rm unor}}$, the 
objects are pairs $(\hatX,\sigma)$ consisting of a complex curve \hatX\ and
an anticonformal involution $\sigma$ that implements 
the action of the Galois group of $\C/\R$. In the case 
of $\calg^{{\rm or}}$, we fix a global section of $\pi$ as an additional datum.

Next, since we are interested in correlation functions depending on insertion 
points, it is natural to consider simultaneously the family $\calm_{g,m}$ of 
all complex curves with marked points that have the same topological type (i.e.,
genus $g$ and number $m$ of marked points) as \hatX. It is convenient to treat 
the positions of the insertion points and the moduli of the complex structure 
on the same footing. The curves that admit an involution $\sigma$ of the same 
type as \hatX\ and for which the marked points are related by $\sigma$
form a submanifold $\calm_{g,m}^\sigma$ of $\calm_{g,m}$.
(To be precise, one obtains \cite{buse} such a relation for Teichm\"uller
spaces, rather than for moduli spaces.)

Given the modular tensor category $\calc$, the complex modular functor 
\cite{BAki} provides us with a vector bundle $\calv$ with projectively
flat connection on $\calm_{g,m}$. 
We can now formulate the `principle of holomorphic factorization' (which for 
certain classes of conformal field theories follows from chiral Ward identities 
that can formally be derived from an action functional \cite{witt39}). 
It states that, first of all, the conformal surface \X\ should be decorated 
in such a way that the double \hatX\ has the structure of an object in the 
decorated cobordism category for the topological field theory based on $\calc$.
It then makes sense to require, secondly, that the correlation function is a 
certain global section of the restriction of $\calv$ to $\calm_{g,m}^\sigma$.

At this point, it proves to be convenient to use the equivalence of the complex
modular functor and the topological modular functor $\tft_\calc$ based on the
modular tensor category $\calc$ \cite{BAki} so as to work in a topological 
(rather than complex-analytic)
category. We are thereby lead to the description of a correlation function on 
\X\ as a specific vector $\Cf(\X)$ in the vector space $\tft_\calc(\hatX)$ that 
is assigned to the double \hatX\ by the topological modular functor $\tft_\calc$.
These vectors must obey two additional axioms:
   \def\leftmargini{1.1em}~\\[-1.2em]
   \begin{itemize} 
\item {\em Covariance\/}: Given any morphism $f{:}\ \X\,{\to}\,\Y$ in the 
      relevant decorated geometric category $\calg_\calc$, we demand
$$ \Cf(\Y) = \tft_\calc(f) \big( \Cf(\X) \big) \, . $$
\item {\em Factorization\/}: Certain factorization properties must be fulfilled.

We refer to \cite{fjfrs,fjfrs2} for a precise formulation of these constraints.
\end{itemize}
The covariance axiom implies in particular that the vector $\Cf(\X)$
is invariant under the action of the mapping class group $\,{\rm Map}(\X)
\,{\cong}\,{\rm Map}(\hatX)^\sigma_{\phantom|}$. This group, also called the 
relative modular group \cite{bisa2}, acts genuinely on $\tft_\calc(\hatX)$.

\section{Surface holonomy}\label{TFTRCFT}

To find solutions to the covariance and factorization constraints on the vectors 
$\Cf(\X)\iN\tft_\calc(\hatX)$ we use the three-dimensional topological field 
theory associated to the modular tensor category $\calc$. Thus we look for 
a (decorated) co\-bor\-dism $(\M_\X,\emptyset, \hat \X)$ such that the vector 
$\tft_\calc(\M_\X,\emptyset,\hatX)1 $\linebreak[0]$ {\in}\, \tft_\calc(\hatX)$
is the correlator $\Cf(\X)$. 

The three-manifold $\M_\X$ should better not introduce any topological 
information that is not already contained in \X. This leads to the idea to use 
an interval bundle as a ``fattening'' of the world sheet.
It turns out that the following quotient of the interval bundle on $\hatX$,
called the connecting (three-)\,manifold, is appropriate \cite{fffs3}:
$$  
\M_\X = \big(\,\hatX \,{\times}\, [-1,1]\,\big) \,{\slash}\,
\langle(\sigma,t\,{\mapsto}\, {-}t)\rangle  \,. 
$$
This three-manifold is oriented, has boundary $\partial\M_\X\,{\cong}\,\hat \X$,
and it contains \X\ as a retract: the embedding $\io$ of \X\ is to the fiber
$t{=}0$, the retracting map contracts along the intervals.

The connecting manifold $\M_\X$ must now be decorated with the help of the 
decoration categories $\HOM(A,A')$. We will describe this procedure for the 
oriented case only. The conformal surface \X\ is decomposed by defect lines 
(which are allowed to end on $\partial\X$) into various two-dimensional regions.
There are two types of one-dimensional structures: boundary components of \X\ 
and defect lines. Defect lines, in general, form a network; they
can be closed or have end points, and in the latter case they can end either
on the boundary or in the interior of \X. Both one-dimensional structures are 
partitioned into segments by marked ``insertion'' points. The end points of 
defect lines carry insertions, too.  Finally, we also allow for insertion 
points in the interior of two-dimensional regions.

To these geometric structures, data of the 2-category $\FROB_\calc$ are now 
assigned as follows. First, we attach to each two-dimensional region a symmetric
special Frobenius algebra, i.e.\ an object of $\FROB_\calc$. To a segment of a 
defect line that separates regions with label $A$ and $A'$, respectively, we 
associate a 1-morphism in $\HOM(A,A')$, i.e.\ an $A$-$A'$-bimodule. Similarly, 
to a boundary segment adjacent to a region labeled by $A$, we assign an object 
in $\HOM(\Cardy,A)$, i.e.\ a left $A$-module. Finally, zero-dimensional 
geometric objects are labeled with 2-morphisms of $\FROB_\calc$; 
in particular, junctions of defect lines with each other or with a boundary 
segment are labeled by 2-morphisms from the obvious spaces.

Two types of points, however, still deserve more comments: those separating 
boundary segments on the one hand, and those separating or creating segments of 
defect lines or appearing in the interior of two-dimensional regions on the 
other. These are the {\em insertion points\/} that were mentioned above.
An insertion point $p\iN\partial\X$ that separates two boundary segments labeled
by objects $M_1,M_2\iN \HOM(\Cardy,A)$ has a single preimage under the canonical 
projection $\pi$ from \hatX\ to \X; to the interval in $\M_\X$ that joins this 
preimage to the image $\io(p)$ of $p$ under the embedding $\io$ of \X\ into 
$\M_\X$, we assign an object $U$ of the category $\calc\eq\HOM(\Cardy,\Cardy)$ 
of chiral data. To the insertion point itself, we then attach a 2-morphism in 
the morphism space $\Hom(M_1\oti U,M_2)$ in $\HOM(\Cardy,A)$.

An insertion point in the interior of \X\ has two preimages on \hatX; 
these two points are connected to $\io(p)$ by two intervals. 
To each of these two intervals we assign an object $U$ and $V$, 
respectively, of the category $\calc$ of chiral data. In the oriented case, the 
global section of $\pi$ is used to attribute the two objects $U,V$ to the two 
preimages. (For the unoriented case, the situation is more involved; in 
particular, the Jandl structure on the relevant Frobenius algebra enters the 
prescription.) We now first consider an insertion point separating a segment of 
a defect line labeled by an object $B_1\iN\HOM(A,A')$ from a segment labeled by 
$B_2\iN\HOM(A,A')$. We then use the left action $\rho_l$ of $A$ and the right 
action $\rho_r$ of $A'$ on the bimodule $B_1$ to define a bimodule structure 
on the object $U\oti B_1 \oti V$ of $\calc$ by taking the morphisms
$(\id_U\oti\rho_l\oti\id_V)\cir (c^{-1}_{U,A}\oti\id_{B_1}\oti\id_V)$ and
$(\id_U\oti\rho_r\oti\id_V)\cir (\id_U\oti\id_{B_1}\oti c^{-1}_{\!A',V})$ as
the action of $A$ and $A'$, respectively, where $c$ denotes the braiding 
isomorphisms of $\calc$. 
The insertion point separating the defect lines is now labeled by a 2-morphism 
in $\Hom(U\oti B_1 \oti V,B_2)$, i.e.\ by a morphism of $A$-$A'$-bi\-mo\-dules.

To deal with insertion points in the interior of a two-dimensional region 
labeled by a Frobenius algebra $A$, we need to invoke one further idea: such a 
region has to be endowed with (the dual of) a triangulation $\Gamma$. 
To each edge of $\Gamma$ we attach the morphism 
$\Delta\cir\eta \iN \Hom(\one,A\oti A)$,
and to each trivalent vertex of $\Gamma$ the morphism 
$\eps \cir m \cir(m \oti \id_A) \iN \Hom(A\oti A\oti A,\one)$. 
\void{
To each 
edge of $\Gamma$ we attach the canonical morphism that identifies the symmetric 
Frobenius algebra $A$ with its dual $A^\vee$, and to each trivalent vertex of 
$\Gamma$ the morphism in $\Hom(A\oti A\oti A,\one)$ that is obtained by 
composing the multiplication morphism of $A$ with the canonical morphism 
$A\,{\to}\,A^\vee$ and a duality morphism. 
}
This pattern is characteristic for 
notions of surface holonomy. It has appeared in lattice topological field 
theories \cite{fuhk} and shows up in the surface holonomy of bundle gerbes as 
well. (For more details, references, and the relation to the Wess-Zumino term 
of WZW conformal field theories in a Lagrangian description, see \cite{gaRe}.) 

Now each of the insertion points $p$ that we still need to discuss is located 
inside a two-dimensional region labeled by some Frobenius algebra $A$ or creates
a defect line. For the first type of points, we choose the triangulation such 
that an $A$-ribbon passes through $p$; to $p$ we then attach a bimodule morphism
in $\Hom(U\oti A \oti V,A)$, with $U$ and $V$ objects of $\calc$ as above.
To a point $p$ at
which a defect line of type $B$ starts or ends, we attach a bimodule morphism 
in $\Hom(U\oti A \oti V,B)$ and in $\Hom(U\oti B \oti V,A)$, respectively.

We have now obtained a complete labelling of a ribbon graph in the connecting
manifold $\M_\X$ with objects and morphisms of the modular tensor category
$\calc$; in other words, a cobordism from $\emptyset$ to \hatX\ in the decorated 
geometric category $\calg_\calc$. Applying the modular functor for the tensor 
category $\calc$ to this cobordism, we obtain a vector
$$ 
\Cf(\X) = \tft_\calc(\M_\X)\,1 \,\in \tft_\calc(\hatX) \,. 
$$
This is the prescription for RCFT correlation functions in the TFT approach.
It follows from the defining properties of a symmetric
special Frobenius algebra that $\Cf(\X)$ does not depend on the choice of
triangulation; for details see \cite{fjfrs}.

\section{Results}

On the basis of this construction one can establish many further results.
Let us list some of them, without indicating any of their proofs:
   \def\leftmargini{4.4em}~\\[-1.2em]
   \begin{itemize}\addtolength\itemsep{2pt}
\item[\cite{fuRs4}] 
     Of particular interest are the correlators for \X\ being the torus
     or the annulus without field insertions, but possibly with defect lines. 
     {}From these ``one-loop amplitudes'' one can derive concrete 
     expressions for partition functions of boundary, bulk and defect fields. 
\\[3pt] The coefficients of these partition functions in the distinguished basis 
     of the zero-point blocks on the torus that is given by characters can be
     shown to be equal to the dimensions of certain spaces of 2-morphisms of the
     2-category $\FROB_\calc$. Thus in particular they are non-negative integers.
\\[3pt] In fact, one recovers expressions that had also been obtained in an
     approach based on subfactors \cite{lore,boek3}. Moreover, these
     coefficients can be shown to satisfy other consistency requirements
     like forming so-called NIMreps of the fusion rules.
\item[\cite{fuRs8}] 
     To extend these results to unoriented (in particular, to unorientable)
     surfaces one must specify as additional datum a Jandl structure on the
     relevant Frobenius algebra. One can then e.g.\ compute the partition
     functions for the M\"obius strip and Klein bottle. Their coefficients in 
     the distinguished basis of zero-point torus blocks are integers, and
     for CFT models which serve as building blocks of type I string theories,
     these partition functions combine with the torus and annulus amplitudes in
     a  way consistent with an interpretation in terms of state spaces of the
     string  theory.
\item[\cite{fuRs9}] 
     The expressions for correlation functions can be made particularly 
     explicit for \cfts\ of simple current type \cite{scya6}, which correspond
     to Frobenius algebras for which every simple subobject is invertible.
     Eilenberg-Mac Lane's \cite{eima2} abelian group cohomology 
     turns out to provide a crucial tool for analyzing this case.
\\[3pt] It should be stressed, though, that the TFT approach to RCFT correlators
     treats the simple current case and other \cfts\ (i.e.\ those having an
     `exceptional modular invariant' as their torus partition function)
     on an equal footing.
\item[\cite{fuRs10}] 
     By expressing some specific correlation functions for the sphere, the disk,
     and the real projective plane through the appropriate (two- or three-point)
     conformal blocks, one can derive explicit expressions for the coefficients
     of operator product expansions of bulk, boundary, and defect fields.
\item[\cite{fjfrs}] 
     For arbitrary topology of the surface \X\ the correlators obtained in the
     TFT construction can be shown to satisfy the covariance and factorization
     axioms that were stated at the end of Section \ref{CORF}.
\item[\cite{ffrs3}] 
     The Picard group of the tensor category $\HOM(A,A)$ describes symmetries
     of the full conformal field theory that is associated to $A$. The fusion
     ring $K_0(\HOM(A,A))$ of that category contains
     information about Kramers\hy Wan\-nier\,-like dualities as well.
\end{itemize}

\section{Conclusions}

The TFT approach to the construction of CFT correlation functions, which 
represents CFT quantities as invariants of knots and links in three-manifolds, 
relates a general paradigm of quantum field theory to the theory of (symmetric 
special) Frobenius algebras in (modular) tensor categories. It thereby 
constitutes a powerful algebraization of many questions that arise in the 
study of conformal field theory. As a result, one can both make rigorous 
statements about rational conformal field theories and set up efficient 
algorithms for the computation of observable CFT quantities.

\smallskip

A rich dictionary relating algebraic concepts and physical notions is emerging. 
It includes in particular the following entries:
   \def\leftmargini{1.1em} \begin{itemize}
\item 
The classification of (oriented) full conformal field theories for given
chiral data $\calc$ amounts to the classification of Morita classes of
Frobenius algebras in $\calc$. \\
As a special case, the classification of those theories whose torus partition
function is ``of automorphism type'' amounts to determining the
Brauer group of the category $\calc$. 
\item 
The Picard group of the tensor category $\HOM(A,A)$ acts as a symmetry group on 
the full conformal field theory associated to $A$, while the fusion ring of 
this tensor category contains information about Kramers-Wannier like dualities.
\item 
Deformations of the full conformal field theory are controlled
by the Davy\-dov-Yet\-ter cohomology of the pair $(\calc,\CA)$.    
\end{itemize}
The structure of this dictionary gives us confidence that some of the insights
of the TFT approach -- though, unfortunately, not most of the proofs -- will still 
be relevant for the study of conformal field theories that are not rational any 
more.

\newpage 

\newcommand\BiA[3]    {\bibitem[#2]{#1} #3}
\newcommand\BiP[3]    {\bibitem{#1} #3}

 \newcommand\wb{\,\linebreak[0]} \def\wB {$\,$\wb}
 \newcommand\JO[6]    {{\em #6}, {#1} {#2} ({#3}), {#4--#5} }
 \newcommand\J[7]     {{\em #7}, {#1} {#2} ({#3}), {#4--#5}\, {{\tt [#6]}}}
 \newcommand\JJ[6]    {{\em #6}, {#1} {#2} ({#3}), {#4}\, {{\tt [#5]}}}
 \newcommand\PhD[3]   {{\em #3}, Ph.D.\ thesis, #1\, {\tt [#2]}}
 \newcommand\Pret[2]  {{\em #2}, pre\-print {\tt #1}}
 \newcommand\BOOK[4]  {{\em #1\/} ({#2}, {#3} {#4})}
 \newcommand\inBO[8]{{\em #8}, in:\ {\em #1}, {#2}\ ({#3}, {#4} {#5}), p.\ {#6--#7}}
 \newcommand\inBo[8]{{\em #8}\ in:\ {\em #1}, {#2}\ ({#3}, {#4} {#5}), p.\ {#6--#7}}
 \def\dim   {dimension}
 \def\jf    {J.\ Fuchs}
 \def\adma  {Adv.\wb Math.}
 \def\anma  {Ann.\wb Math.}
 \def\coma  {Con\-temp.\wb Math.}
 \def\comp  {Com\-mun.\wb Math.\wb Phys.}
 \def\cOmp  {Com\-mun. Math.\wb Phys.}
 \def\Comp  {Com\-mun.\wb Math. Phys.}
 \def\cpma  {Com\-pos.\wb Math.}
 \def\fiic  {Fields\wB Institute\wB Commun.}
 \def\ihes  {Publ.\wb Math.\wB I.H.E.S.}
 \def\ijmp  {Int.\wb J.\wb Mod.\wb Phys.\ A}
 \def\jhep  {J.\wb High\wB Energy\wB Phys.}
 \def\joal  {J.\wB Al\-ge\-bra}
 \def\jomp  {J.\wb Math.\wb Phys.}
 \def\josp  {J.\wb Stat.\wb Phys.}
 \def\jpaa  {J.\wB Pure\wB Appl.\wb Alg.}
 \def\maze  {Math.\wb Zeitschr.}
 \newcommand\ncmp[3] {\inBO{IXth International Congress on
            Mathematical Physics} {B.\ Simon, A.\ Truman, and I.M.\ Da\-vis,
            eds.} {Adam Hilger}{Bristol}{1989} {#1}{#2}{#3} }
 \def\nuci  {Nuovo\wB Cim.}
 \def\nupb  {Nucl.\wb Phys.\ B}
 \def\Nupb  {Nucl. Phys.\ B}
 \def\phlb  {Phys.\wb Lett.\ B}
 \def\phrl  {Phys.\wb Rev.\wb Lett.}
 \def\pnas  {Proc.\wb Natl.\wb Acad.\wb Sci.\wb USA}
 \def\pspm  {Proc.\wb Symp.\wB Pure\wB Math.}
 \def\ptrs  {Phil.\wb Trans.\wb Roy.\wb Soc.\wB Lon\-don}
 \def\rims  {Publ.\wB RIMS}
 \def\rvmp  {Rev.\wb Math.\wb Phys.}
 \def\taia  {Topology\wB and\wB its\wB Appl.}
 \newcommand\tgqf [3] {\inBO{Topology, Geometry and Quantum Field Theory
            {\rm [London Math.\ Soc.\ Lecture Note Series \#~308]}}
            {U.\ Tillmann, ed.} \CUP\Ca{2004} {#1}{#2}{#3}}
 \newcommand\tgqF [3] {\inBo{Topology, Geometry and Quantum Field Theory
            {\rm [London Math.\ Soc.\ Lecture Note Series \#~308]}}
            {U.\ Tillmann, ed.} \CUP\Ca{2004} {#1}{#2}{#3}}
 \newcommand\ndha[3] {\inBO{New Directions in Hopf Algebras {\rm [MSRI
            Publications No.\ 43]}} {S.\ Montgomery and H.-J.\ Schneider, eds.}
            \SV\NY{2002} {#1}{#2}{#3}}
 \def\trgr  {Trans\-form. Groups}
   \def\AMS    {{American Mathematical Society}}
   \def\BIR    {{Birk\-h\"au\-ser}}
   \def\Bo     {{Boston}}
   \def\Ca     {{Cambridge}}
   \def\CUP    {{Cambridge University Press}}
   \def\NY     {{New York}}
   \def\PR     {{Providence}}
   \def\SV     {{Sprin\-ger Ver\-lag}}

\end{document}